\documentclass{svmult}

\usepackage{makeidx} 
\usepackage{graphicx}
\usepackage{multicol}
\usepackage[bottom]{footmisc}
\usepackage{natbib}
\usepackage{url}
\usepackage{subfigure}
\renewcommand{\PackageWarningNoLine}[2]{}
\usepackage{amsmath}
\usepackage{amsfonts}
\usepackage{wrapfig}

\newcommand{\amat}{\arrfont{C}}                 
\newcommand{\arrfont}[1]{\mbox{\sffamily$\textbf{#1}$}}
\newcommand\assign{\mathbin{:=}}                
\newcommand{\be}{\begin{equation}}
\newcommand{\cpress}{p}                         
\newcommand{\Dmat}{\arrfont{D}}
\newcommand{\DTmat}{\arrfont{D}^{\rm T}}
\newcommand{\ee}{\end{equation}}
\newcommand{\eq}[1]{(\ref{#1})}

\newcommand{\fig}[1]{figure\ \ref{#1}}

\newcommand{\gnum}{\lavec{g}}                   
\newcommand{\Hone}{\set{H}^1}                   

\newcommand{\lavec}[1]{{\boldsymbol#1}}         
\newcommand{\Lmat}{\arrfont{L}}                 
\newcommand{\Ltwo}{\set{L}_2}                   
\newcommand{\Mmat}{\arrfont{M}}                 
\newcommand{\PD}[2]{\partial_{#2}#1}            
\newcommand{\pdomain}{\set{D}}                  
\newcommand{\pnum}{\hat{\lavec{p}}}                   
\newcommand{\polynomialsset}[1]{\set{P}_{#1}}   
\newcommand{\pspace}{\set{Y}}                   

\newcommand{\set}[1]{\mathbf{#1}}               
\newcommand{\setdef}[2]{\left\{#1\left\bracevert#2\right.\right\}}
\newcommand{\si}{\mu}                           

\newcommand{\uspace}{\set{U}}                   
\newcommand{\uv}[1]{{\vec{e}}^{\:#1}}           
\newcommand{\znum}{\hat{\lavec{Z}}}             

 \def\B{{\bf B}} \def\j{{bf j}} \def\b{{\bf b}} \def\u{{\bf u}}   \def\Z{{\bf Z}} \def\w{{\bf w}} \def\ztst{{\bf \zeta}}

\begin{document}

\title*{Optimized Schwarz preconditioning for SEM based magnetohydrodynamics}
\titlerunning{Optimized Schwarz preconditioning for SEM based MHD}

\author{Amik St-Cyr\inst{1}\and
Duane Rosenberg\inst{1}\and Sang Dong Kim\inst{2}}
\institute{
Institute for Mathematics Applied to Geosciences, National Center for Atmospheric 
Research, Boulder, CO, USA: \texttt{\{amik,duaner\}@ucar.edu}\and 
Department of Mathematics, Kyungpook National University, Daegu 702-701, 
South Korea: \texttt{skim@knu.ac.kr}}
%
%
\maketitle

\noindent {\bf Summary.} A recent theoretical result on optimized Schwarz algorithms 
demonstrated at the algebraic level enables the modification of an existing Schwarz procedure to 
its optimized counterpart. In this work, it is shown how to modify a bilinear FEM based Schwarz 
preconditioning strategy originally presented in \cite{Fischer:1997ys} to its optimized version. 
The latter is employed to precondition the pseudo--Laplacian operator arising from the 
spectral element discretization of the magnetohydrodynamic equations in  Els\"asser form.

\section{Introduction}
This work concerns the preconditioning of a pseudo--Laplacian 
operator\footnote{A.k.a: consistent Laplacian or approximate pressure Schur complement.} 
associated with the saddle point problem arising at each time-step in a spectral element 
based adaptive  MHD solver.  The approach proposed herein is a modification of the 
method developed in \cite{Fischer:1997ys} where an overlapping Schwarz preconditioner was 
constructed using a low oder discretization. The latter approach is based on the spectral 
equivalence between finite-elements and spectral elements \cite{Canuto:2007zi,Kim:2006kl}. 
The finite-element blocks, representing the additive Schwarz, are replaced by so called 
optimized Schwarz blocks \cite{St-Cyr:2007dq}. Two types of meshes, employed to 
construct the $\mathbb{Q}_1$ block preconditioning are investigated. The first one is cross 
shaped and shows good behavior for additive Schwarz (AS) and restricted additive Schwarz 
(RAS). However improved convergence rates of the optimized RAS (ORAS) version are 
completely dominated by the corner effects \cite{Chniti:pr}. Opting for a 
second grid that includes the corners seems to correct this issue. For the zeroth order 
optimized transmission condition (O0) an exact tensor product form is available while for 
the O2 version a slight error is introduced in order to preserve the properties of the 
operators and enable the use of fast diagonalization techniques introduced for 
spectral elements in \cite{Couzy:1995tv}. It is shown how to modify trivially an existing 
fast-diagonalization procedure and numerical experiments demonstrate the 
efficiency of the modification.

\section{Governing equations and discretization}
\label{sec:1}
For an incompressible fluid with constant mass density $\rho_0$, the magnetohydrodynamic  
(MHD) equations are:

\begin{equation}
\partial_t {\u} + {\u}\cdot {{\nabla}} {\u} = - {{ \nabla}} \cpress + \j \times \b
    + \nu \nabla^2 {\u} ,
\label{eq_momentum} \end{equation}
\begin{equation}
\partial_t {\b}  = { { \nabla}} \times (\u \times \b)
    + \eta \nabla^2 {\b} 
\label{eq_induction} \end{equation}
\begin{equation}
{ {\nabla}} \cdot {\u} =0, \hskip0.2truein { {\nabla}} \cdot {\b} =0
\label{eq_incompressible}
\end{equation}
where ${\u}$ and $\b$ are the velocity and magnetic field (in Alfv\'en velocity units,
$\b=\B/\sqrt{\mu_0 \rho_0}$ with $\B$ the induction and $\mu_0$ the permeability); 
$\cpress$ is the pressure divided by the mass density, and $\nu$ and $\eta$ are the 
kinematic viscosity and the magnetic resistivity. In  Els\"asser form, the equations are 
\cite{Elsasser:1950vn}:
\be
\partial_t \Z^{\pm} + \Z^{\mp}\cdot { {\nabla}} \Z^{\pm} + { {\nabla}} \cpress -\nu^{\pm}\nabla^2 \Z^{\pm} 
- \nu^{\mp} \nabla^2 \Z^{\mp} = 0
\label{eq_zmomentum}
\ee
\be
 { {\nabla}} \cdot  \Z^{\pm} = 0\,,
\label{eq_zconstraint}
\ee
with $ \Z^{\pm} = \u \pm \b $ and $ \nu^\pm = \frac{1}{2}(\nu \pm \eta)$. The velocity $\u$ and 
magnetic field $\b$ can be recovered by expressing them in terms of $\Z^{\pm}$. Since all 
time-scales are of interest, an explicit second order Runge-Kutta scheme is applied to discretize the 
time-derivative of the above system while, for the spatial part, a 
$\polynomialsset{N}-\polynomialsset{N-2}$  spectral element formulation  was chosen to prevent the 
excitation of spurious pressure modes using the spaces 
\begin{eqnarray}
\uspace_{{\bf \gamma }}&\assign%
\setdef{\w=\sum_{\si=1}^dw^\si\uv{\si}\,}{w^\si\in\Hone({D})\; \forall \mu\;\&\;\w={\bf \gamma}\;{\rm on}\;\partial D}
\end{eqnarray}
\begin{eqnarray}
\Hone(D)&\assign\setdef{f\,}{f\in\Ltwo(D)\;\&\;\PD{f}{x^\si}\in\Ltwo(D)\;\forall\mu} \, ,
\end{eqnarray}
with $\w = \u,\, \b$, $\Z^\pm$, $\cpress$ and their test functions, $\ztst^\pm$ and $q$ restricted to 
finite--dimensional subspaces of these spaces:
\begin{align}
&\Z^\pm    \in \uspace^{N}    \;\;\;= \uspace_{\Z_0}    \; \bigcap \; \polynomialsset{N},&\;
&\ztst^\pm \in \uspace_{0}^{N}\;\;\;\;= \uspace_{{\bf 0}} \; \bigcap \; \polynomialsset{N},&\notag\\
&\cpress, q      \in \pspace^{N-2}   = \Ltwo(D)    \; \bigcap \; \polynomialsset{N-2}& \notag
\end{align}
see for instance  \cite{Y.-Maday:1992zr, Fischer:1997ys}.
The basis for the velocity expansion in $\polynomialsset{N}$ is the set of
Lagrange interpolating polynomials on the Gauss-Lobatto-Legendre ({\rm GL}) quadrature 
nodes, and the basis for the pressure is the set of Lagrange interpolants on the Gauss-Legendre 
({\rm G}) quadrature nodes. In the spectral element method formalism, the domain, $\pdomain$, is 
composed of a union of non-overlapping subdomains, or elements, $\set{E}_k$: 
$ \pdomain = \bigcup_{k=1}^{K} \set{E}_k$, and functions in $\uspace^{N}$ and 
$\pspace^{N-2}$  are represented as expansions in terms of tensor products of basis functions 
within each subdomain. The complete discretization at each stage is:
\be
\znum_{j}^\pm = \znum_{j}^{\pm,n}  - \frac{1}{k} \Delta t \; \Mmat^{-1}(\Mmat\amat^{\mp} \znum_j^\pm  
-  \DTmat_j \pnum^\pm +\nu_{\pm}\Lmat\znum_j^\pm +\nu_{\mp}\Lmat\znum_j^\mp ) .
\label{eq_rkstage}
\ee
We require that each Runge-Kutta stage obeys \eq{eq_zconstraint} in its discrete form, so multiplying
\eq{eq_rkstage} by $\Dmat_j$, summing, and setting the term $\Dmat^j \znum_j^\pm = 0$
leads to the
following pseudo--Poisson equation for the pressures, $\pnum^\pm$:
\be
\Dmat^j \Mmat^{-1} \DTmat_j \pnum^\pm = \Dmat^j \gnum_j^\pm \, ,
\label{eq_pseudopoisson}
\ee
where the quantity
$$
\gnum_j^\pm =  \frac{1}{k} \Delta t \; \Mmat^{-1}(\Mmat \amat^{\mp}\znum_j^\pm  + 
\nu_{\pm}\Lmat\znum_j^\pm +\nu_{\mp}\Lmat\znum_j^\mp ) - \znum_j^{\pm,n}.   
$$
is the remaining inhomogeneous contribution (see \cite{D.-Rosenberg:2007bh}). More 
details on the various operators can be found in \cite{Deville:2002kx}. Equation 
\eq{eq_pseudopoisson} is solved using a preconditioned 
iterative Krylov method.

\section{From classical to optimized Schwarz}
\label{sec:2}
 The principle behind optimized Schwarz methods consists into replacing 
the Dirichlet {\em transmission} condition present in the classical Schwarz approach 
by a more general boundary condition. This 
idea was first analyzed by Lions in \cite{Lions:1988:SAM} where a Robin condition was 
introduced. The latter contained a positive parameter that could possibly be used to enhance 
convergence. However, until recently, it was not clear how to define optimally that parameter 
for the new conditions at the interfaces between subdomains. Optimized Schwarz methods 
are derived from a Fourier analysis of the continuous elliptic partial differential equation, 
see for instance \cite{gander:699} and references therein. Starting from the problem stated at the continuous level, 
suppose that a linear elliptic operator $\mathcal{L}$ with forcing $f$ and boundary conditions 
$\mathcal{P}$ needs to be solved on $\pdomain$. An algorithm that can be employed to solve 
the global problem $\mathcal{L} u = f$ is 
\begin{align}
\label{JacobiSchwarz}
  \mathcal{L}u^{n+1}_i & =  f \;\;\mbox{in $ \set{\bar{E}}_i$} \notag\\
  \mathcal{P}(u^{n+1}_i) & =  g \;\;\mbox{on $\partial \pdomain \cap \set{\bar{E}}_i$}\\
  u^{n+1}_i & =  u^{n}_j \;\;\mbox{on $\Gamma_{ij}$}\notag
 \end{align}
where the sequence with respect to $n$ will be convergent for any initial guess $u^0$ with 
$$\Gamma_{ij}\equiv \partial\set{\bar{E}}_i \cap \set{\bar{E}}_j \ne\{ \emptyset\}.$$
\begin{wrapfigure}{l}{0.62\textwidth}
  \includegraphics[width=0.65\textwidth]{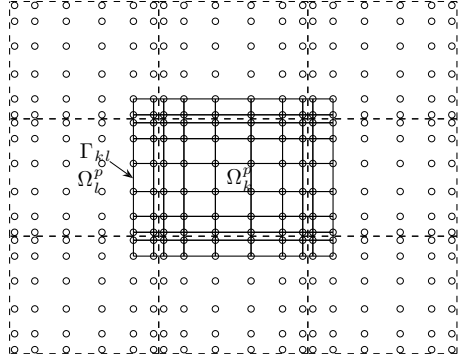}
  \caption{One overlapping element}\label{fig:ovlelem}
\end{wrapfigure}

This is none other than the classical Schwarz algorithm at the continuous level corresponding 
to RAS at the matrix level. In Fig. \ref{fig:ovlelem}, $\Omega^p_k$ 
represents the quadrangulation of the overlapping domain $\set{\bar{E}}_k$.  
The optimized version of the above algorithm replaces the transmission conditions between subdomains by 
\begin{equation}\label{JacobiSchwarzOpt}
  \begin{array}{rcll}
  \tilde{\mathcal{B}}_{ij}u^{n+1}_i & = &\tilde{\mathcal{B}}_{ij}u^{n}_j ;\;\mbox{on $\Gamma_{ij}$}\notag\\
  \end{array}
\end{equation}
where $\tilde{\mathcal{B}}_{ij}$ makes each subproblem well posed and can be a function of optimizable parameters: the algorithm, like in the classical case, converges to the solution 
of $ \mathcal{L}u$ with $\mathcal{P}(u) = g$ on $\pdomain$. The discrete algebraic 
version is   
\[
\left(
\begin{array}{cc}
 A^k_{ii} & A^k_{i \Gamma}     \\
 C^k_{\Gamma i} & C^{k}_{\Gamma\Gamma}       
\end{array}
\right)
\left(
\begin{array}{c}
 u^k_{i}      \\
 u^k_{\Gamma}         
\end{array}
\right)^{n+1}=
\left(
\begin{array}{c}
 f^k_{i}      \\ 
C^k (u_{\partial\Omega^p_k})^n         
\end{array}
\right)
\]
with $C^k_{\Gamma i}$, $C^{k}_{\Gamma\Gamma}$ and $C$ corresponding to the 
the discrete expressions  of the new transmission conditions. At this point notice that 
$A^k_{ii}$ is exactly the same block as in the original Schwarz algorithm. A simple  
(block) Gaussian elimination leads to the following preconditioned 
system
\begin{align}
\{I  -  \sum^K_{j,k=1} (\tilde{R}^k)^T \tilde{B}_{kj} R^j\} u =  \sum^K_{k=1} (\tilde{R}^k)^T (\tilde{A}^k_{ii})^{-1}R^k f 
\end{align}
with $\tilde{A}^k_{ii} = A^k_{ii} - A^k_{i \Gamma} (C^{k}_{\Gamma\Gamma})^{-1} C^k_{\Gamma i}$.
Thus, the optimized restricted additive Schwarz preconditioning is  expressed as 
$P^{-1}_{ORAS} \equiv  \sum^K_{k=1} (\tilde{R}^k)^T (\tilde{A}^k_{ii})^{-1}R^k $.
The above results are completely algebraic and independent of the underlying 
space discretization method. The complete proof in the additive and multiplicative 
case with and without overlap can be found in \cite{St-Cyr:2007dq}. In the case of two 
subdomains it can be shown that the optimal transmission operator is the Schur complement 
\cite{St-Cyr:2007cr}. Also in weighted residual techniques the artificial boundary conditions
should be properly weighted.

\section{$Q_1$ formulation and optimized Schwarz}
The $Q_1$ formulation for the Schwarz algorithm can be found in (cite).  We merely express 
here the changes necessary in order to obtain an optimized version. First, the overlap depicted 
in Fig. \ref{fig:ovlelem} is the minimal one. Secondly, the normal or tangential derivatives expressions 
at the boundaries must not involve more than $2$ points.  Including more would destroy 
the optimized iterates as mentioned in \cite{St-Cyr:2007dq} (remark $1$).  Both 
requirements are satisfied by the $Q_1$ FEM formulation naturally.

\begin{figure}
  \centering
  \includegraphics[scale=0.5]{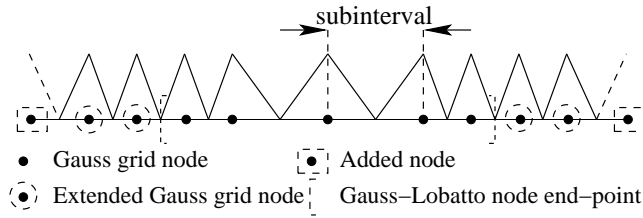}
  \caption{Schematic of Q1 assembly. The left-- and right--most dashed shape functions
           are not included in the assembly.}\label{fig:q1assembly}
\end{figure}

It is important to handle the assembly of the Q1 operators correctly at the endpoints.
Thus we build our reference interval on the extended grid to include
one additional node at each end of the interval. We then begin the assembly (direct stiffness summation)
starting at the second node on the left, and continuing until the second-to-last node on the right. 
In this way, the negative--sloping 
linear FEM shape function on the left--most subinterval and the positive--sloping shape function 
on the rightmost subinterval are not included in the assembly. The general idea is illustrated 
in \fig{fig:q1assembly}.  

We can define for the linear problem in equation \eq{JacobiSchwarz} a general transmission condition \cite{St-Cyr:2007dq} 
between each element
\begin{align*}
\left[ \frac{\partial u_j}{\partial {\bf n}} + T(u_j,p,q,\tau) 
\right]^{n+1}_{\Gamma_{ij}} = 
\left[\frac{\partial u_i}{\partial {\bf n}} +T(u_j,p,q,\tau)
\right]^n_{\Gamma_{ij}}
\label{trnscond}
\end{align*}
where the interface blocks $T(u_j,p,q,\tau)\equiv p u_j - q \frac{\partial^2 u_j}{\partial
  \tau^2}$,
define a transmission condition of order $2$ with two parameters, $p$ and $q$, specified
in \cite{gander:699}, and also provided for completeness in table \ref{tab1}. In general, 
$p$ and $q$ are different depending on whether there is overlap, but here we consider only the case of
finite overlap.
\begin{table}
  \centering
  \begin{tabular}{|l|c|c|}\hline
        &  $p$   &  $q$ \\\hline
   OO0, overlap $Ch$ &
   $2^{-1/3}(k^2_{\min}+\eta)^{1/3}(Ch)^{-1/3}$ & 0 \\
   OO2, overlap $Ch$ &
    $2^{-3/5}(k_{\min}^2+\eta)^{2/5}(Ch)^{-1/5}$ &
    $2^{-1/5}(k_{\min}^2+\eta)^{-1/5}(Ch)^{3/5}$ \\\hline
  \end{tabular}
  \caption{Choices for the parameters $p$ and $q$ used in the interface
blocks.  OOj stands for optimized of order j. }
  \label{tab1}
\end{table}

\subsection{FDM}
When rectangular elements are considered a {\it fast diagonalization method} (FDM) (e.g., \cite{Couzy:1995tv})
can be used to invert the 
optimized $Q_1$ blocks. The number of operations required to invert $N^d \times N^d$ 
matrix using such a technique is $O(N^{d+1})$ and the application of the inverse is performed 
using efficient tensor products in $O(N^d)$ operations.  We propose the form 
\begin{align}
\tilde{A}_{ii} = (M + T_0 q)\otimes (K + T_0 p) + (K + T_o p) \otimes (M+ T_0 q)
\end{align}
as the optimized block matrix where $T_0$ is a matrix almost completely filled with zeros 
with the exception of the entries $(1,1)$ and $(N,N)$ which are set to $1$. 
This form implements the weak form of the transmission conditions, Equation \eq{trnscond}, directly
in the optimized block, $\tilde{A}_{ii}$. Notice that in order 
for the fast diagonalization technique to 
be applicable, the coefficients $p$ and $q$ must be constant on their face. The modified matrix  
$M+T_o q$ is still symmetric and positive definite while the matrix $K + T_o p$ is still symmetric.
This enables the use of the modified mass matrix in an inner product and the simultaneous 
diagonalization of both tensors. When $q=0$ the proposed formula is exact; however, 
when $q\ne0$
a slight error is made at the corners. If the quantity $(T_0 p)\otimes (T_0 q)$ was removed 
then the expression in the $O2$ would also be exact.

\section{Numerical experiments}
\begin{figure}
\begin{center}
\mbox{ \includegraphics*[width=0.49\textwidth]{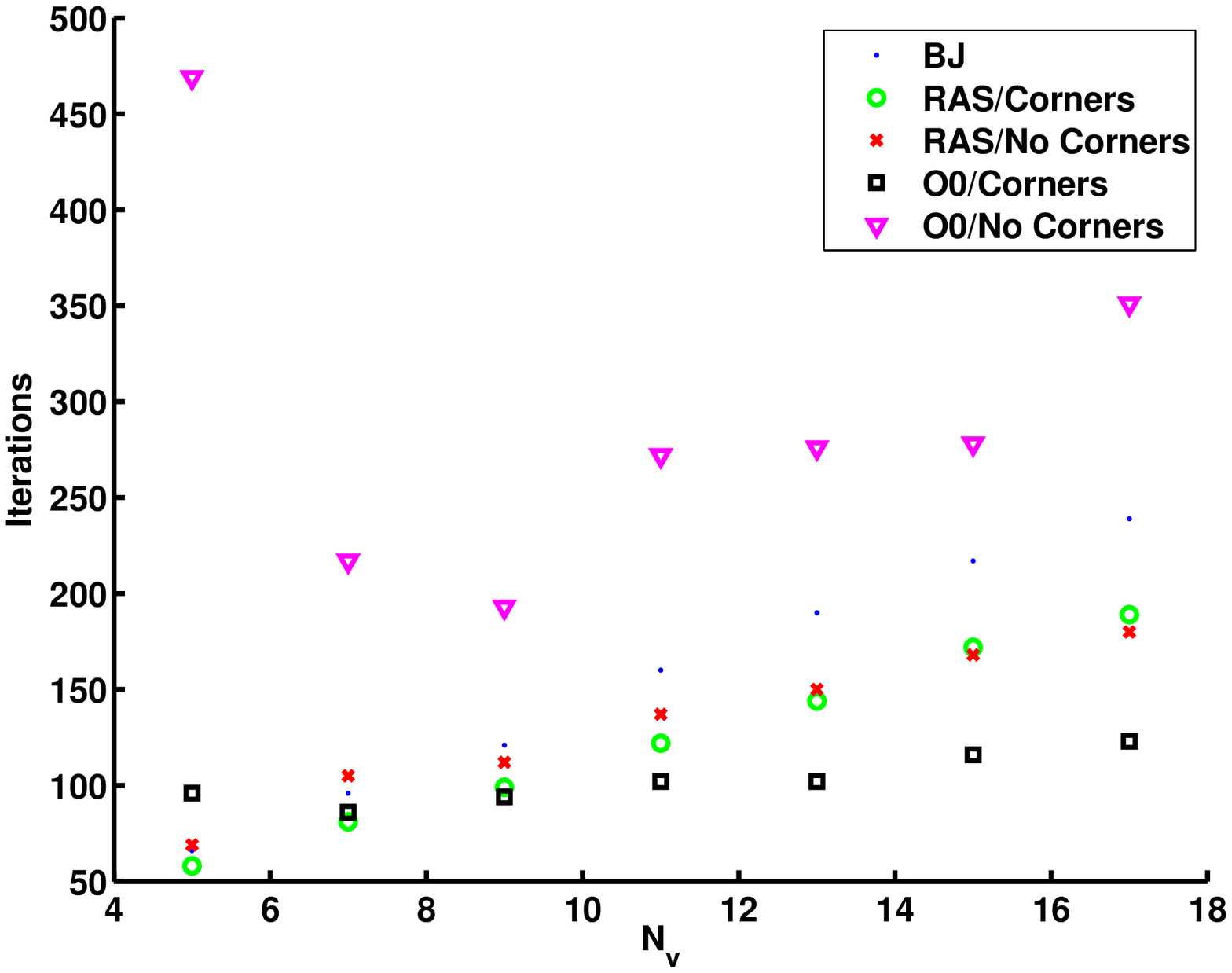}
       \includegraphics*[width=0.49\textwidth]{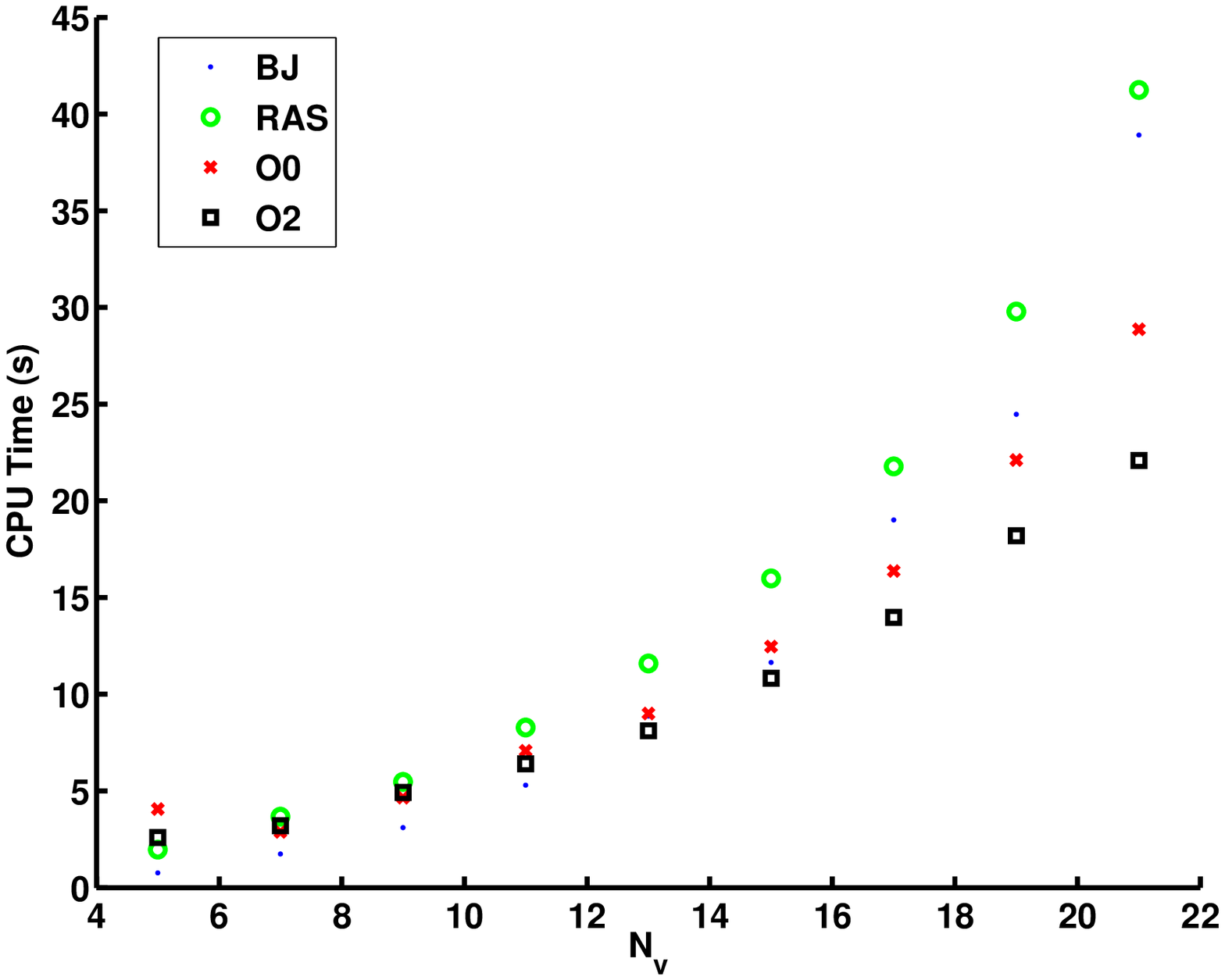} }
       \end{center}
\caption{{\it left:} Plot of iteration count vs GLL-expansion node number for different preconditioners with
and without corner communication on an $8\times8$ element grid.
{\it right:} Comparison of of CPU time vs. GLL-expansion node number of FDM-based preconditioners with 
corner communication on a $16\times16$ element grid.
}\label{fig:cornerperfchk}
\end{figure}

We have implemented the RAS preconditioner described above in the MHD code. This 
version allows for variable overlap of the extended grid. 
The ORAS counterpart has also been implemented as described by starting from the RAS, and
for comparison, we can use a high-order block Jacobi (BJ) method as well.
%
We consider first tests of a single pseudo--Poisson solve. We use a periodic grid of
$E=8\times8$ elements, and iterate using BiCGStab until the residual is $10^{-8}$ times
that of the initial residual. The extended grid overlap is $2$, and the initial starting guess 
for the Krylov method is composed of random noise. 
The first test, uses non-FDM preconditioners to investigate 
the effect of including corner transfers on the optimization. The results are
presented in Fig. \ref{fig:cornerperfchk}, in which we consider only the O0 optimization. Note that
even thought the RAS is much less sensitive to the corner communication, especially at higher
$N_v$, the O0 with corners requires fewer 
iterations. Clearly, corner communication is crucial to the proper functioning of the optimized methods.
%

In the next experiment, all the parameters are maintained except we use a grid of $E=16\times16$ elements
together with the FDM version of the preconditioners to 
investigate performance. These results are presented on the right-most figure of Fig. \ref{fig:cornerperfchk}.

\section{Conclusions and future directions}
It is shown that a simple modification of a RAS in a low order FDM based 
preconditioning of the the pseudo--Laplacian operator can reduce the time 
to solution by up to a factor of two for high order GLL expansions. Also, as expected 
from the work of \cite{Chniti:pr}, we find that the cross form of the subdomains is not 
suitable as is for the optimized version 
of the algorithm: corners need to be included. Upcoming work will concern the 
inclusion of a coarse solver in this approach and a treatment for non-conforming
elements.

\subsection*{Acknowledgments}
The third author was upported by Korea's Research Foundation under grant number 
{\tt KRF-2005-070-C00017}. The first author is grateful for the support provided 
by the BK21 program. NCAR is supported by the National Science Foundation.
%
%
\bibliographystyle{alpha}
\bibliography{references}



\end{document}